\documentclass[11pt]{article}
\usepackage{amscd, amsmath, amssymb}
\usepackage[latin10]{inputenc}

\title{Jumping of the nef cone for Fano varieties}
\author{Burt Totaro}
\date{  }

\def\Z{\text{\bf Z}}
\def\Q{\text{\bf Q}}
\def\R{\text{\bf R}}
\def\C{\text{\bf C}}
\def\P{\text{\bf P}}

\def\arrow{\rightarrow}

\def\imp{\Rightarrow}

\def\qed{\ QED }

\def\Pic{\text{Pic}}
\def\Hom{\text{Hom}}

\def\dash{\dashrightarrow}
\def\rank{\text{rank}}
\def\tr{\text{tr}}
\def\Cl{\text{Cl}}

\def\sm{\text{sm}}

\def\tx{\widetilde{X}}

\def\A5{\widetilde{A^5}}
\def\T{{\cal T}}
\def\Nef{\text{Nef}}
\def\Cox{\text{Cox}}

\def\Kb{\overline{K}}
\def\Gal{\text{Gal}}
\def\Spec{\text{Spec}}
\def\Zh{\widehat{\Z}}
\def\et{\text{et}}
\def\Hilb{\text{Hilb}}
\def\Proj{\text{Proj}\, }
\def\im{\text{im}}

\begin{document}
\maketitle

\newtheorem{theorem}{Theorem}[section]
\newtheorem{corollary}[theorem]{Corollary}
\newtheorem{lemma}[theorem]{Lemma}

Among all projective algebraic varieties, Fano varieties (those with ample
anticanonical bundle) can be considered the simplest.
Birkar, Cascini, Hacon, and M\textsuperscript{c}Kernan
showed that the Cox ring of a Fano variety,
the ring of all sections of all line bundles, is finitely generated
\cite{BCHM}. This implies a fundamental fact about the birational
geometry 
of a Fano variety: there are only finitely many small $\Q$-factorial
modifications of the variety, parametrized by a chamber
decomposition of the movable
cone into rational polyhedral cones (the nef cones of the modifications).
See Koll\'ar-Mori \cite{KMbook} or Hu-Keel \cite{HK} for definitions.

There are also strong results about deformations of Fano varieties.
For any deformation of a $\Q$-factorial
terminal Fano variety $X_0$, de Fernex and Hacon showed
that the Cox ring deforms in a flat family
\cite[Theorem 1.1, Proposition 6.4]{DH}. This had been proved for smooth Fanos
by Siu \cite[Corollary 1.2]{Siu}.
In other words, all line bundles have the ``same number'' of sections
on $X_0$ as on deformations of $X_0$. It follows that the
movable cone remains constant when $X_0$ is deformed \cite[Theorem 6.8]{DH}.
To be clear, we use ``deformation'' to mean a nearby deformation;
for example, the statement on the movable cone means that
the movable cone is constant for $t$ in some open neighborhood of $0$,
for any flat family $X\arrow T$ with fiber $X_0$ over a point $0\in T$.

When a
$\Q$-factorial terminal Fano variety $X_0$
is deformed, de Fernex and Hacon asked whether
the chamber decomposition of
the movable cone also remains constant \cite[Remark 6.2]{DH}.
This would say in particular
that the nef cone, or dually the cone of curves, remains
constant under deformations of $X_0$.
The answer is positive in dimension at most 3, and also
in dimension 4 when $X_0$ is Gorenstein \cite[Theorem 6.9]{DH}.
In any dimension, Wi\'sniewski
showed that the nef cone of a {\it smooth }Fano variety remains constant
under deformations \cite{Wisniewskidef, Wisniewskirig}. 

In this paper, we show 
that the blow-up of $\P^4$ along a line
degenerates to 
a $\Q$-factorial terminal Fano
4-fold $X_0$ with a strictly smaller nef cone (Theorem \ref{four}).
Therefore the results by de Fernex and Hacon on deformations of 
3-dimensional Fanos are best possible. The example is based on
the existence of high-dimensional flips which deform to isomorphisms,
generalizing the Mukai flop.
This phenomenon will be common,
and we give a family of examples in various dimensions,
including a Gorenstein example in dimension 5 (Theorem 
\ref{jump}). The examples
also disprove the ``volume criterion
for ampleness'' on $\Q$-factorial terminal Fano varieties
\cite[Question 5.5]{DH}.

In view of Wi\'sniewski's theorem, it would be interesting to
describe the largest class of Fano varieties for which the nef
cone remains constant under deformations. It will not be enough
to assume that the variety is $\Q$-factorial and terminal, but
the optimal assumption should also apply to many Fanos which are
not $\Q$-factorial or terminal.

For example, de Fernex
and Hacon proved that $\Q$-factorial terminal toric Fano varieties
are rigid, which implies that the nef cone is constant
under deformations in that case \cite[Theorem 7.1]{DH}.
In this paper, we show more generally
that a toric Fano variety which is smooth in codimension 2 and
$\Q$-factorial in codimension 3 is rigid (Theorem \ref{toric}).
(This terminology means that the variety is smooth outside some
closed subset of codimension at least 3, and $\Q$-factorial outside
some closed subset of codimension at least 4.)

We show in Theorem \ref{class}
that the divisor class group is constant under deformations of klt
Fano varieties which are smooth in codimension 2 and $\Q$-factorial
in codimension 3. This extends results of Koll\'ar-Mori \cite{KMflip}
and de Fernex-Hacon \cite{DH}.

We show that for any deformation of a terminal
Fano variety which is $\Q$-factorial in codimension 3,
the Cox ring (of all Weil divisors) deforms in a flat family
(Theorem \ref{flat}). The proof uses de Fernex-Hacon's extension theorem
on the $\Q$-factorial case. On the other hand, section \ref{counter}
shows that the Cox ring need not vary in a flat family for Fano
varieties with slightly worse singularities, thus answering
a question in the first version of this paper. In a sense, section
\ref{counter} shows that, for Fano varieties,
de Fernex-Hacon's extension theorem is optimal.

A side result which seems to be new
is that for a complex projective
$n$-fold $X$ with rational singularities such that $H^1(X,O)=H^2(X,O)=0$,
the divisor class group maps isomorphically
to the homology group $H_{2n-2}(X,\Z)$
(Theorem \ref{top}). For 3-folds with isolated singularities, this
was proved by Namikawa-Steenbrink \cite[Theorem 3.2]{NS}.

Thanks to Tommaso de Fernex, Christopher Hacon, and J\'anos Koll\'ar
for useful conversations. The excellent referees simplified the proof
of Theorem \ref{class} and spotted a misuse of Ein-Popa's extension
theorem in my earlier proof of Theorem \ref{flat}.

\section{Jumping of the nef cone}

\begin{theorem}
\label{jump}
For any positive integers $a\geq b$, there is a smooth variety $X$
of dimension $a+b+1$
and a flat projective morphism $t:X\arrow A^1$ with the following
properties. The fibers $X_t$ for $t\neq 0$ are isomorphic
to $\P^a\times \P^b$, and the nef cone of $X_0$ is a proper subset 
of the nef cone of $X_t=\P^a\times \P^b$.
 The fiber $X_0$ is a terminal Gorenstein Fano variety for $a>b$
and is $\Q$-factorial if $a>b>1$. For $a=b$, $X_0$ is a smooth
weak Fano ($-K_{X_0}$ is nef and big).
\end{theorem}

In particular, for $a=3$ and $b=2$, $X_0$ is a $\Q$-factorial terminal
Gorenstein Fano
5-fold whose nef cone changes when $X_0$ is deformed
to $X_t=\P^3\times \P^2$. In fact, the only singular point of $X_0$
is a node. This shows the optimality of de Fernex and Hacon's theorem
that the nef cone is constant under deformation for $\Q$-factorial terminal
Gorenstein Fano 4-folds \cite[Theorem 6.9]{DH}.
Theorem \ref{jump}
works over any field.

{\bf Proof. }We first describe why the nef cone changes
in these examples. The smooth variety $X$ has
a projective morphism $X\arrow Z$ over $A^1$
which contracts a copy of $\P^a$ contained in the special fiber $X_0$.
This is a flipping contraction of $X$ for $a>b$ and a flopping contraction
for $a=b$ (in particular, the canonical bundle
$K_X$ has degree negative or zero, respectively,
on any curve in $\P^a$). Let $Y\arrow Z$ be the corresponding flip
or flop; in our example, $Y$ is also smooth, with the $\P^a$ in $X$
replaced by a $\P^b$ in $Y$. Then the birational map
$X\dash Y$ restricts to a birational map $X_0\dash Y_0$
which is itself a flip or flop, since $K_{X_0}=(K_X+X_0)|_{X_0}=(K_X)|_{X_0}$.
(For $b=1$, $X_0\arrow Z_0$ contracts
a divisor $\P^a$ in $X_0$, but this divisor is not $\Q$-Cartier
and $X_0\dash Y_0$ must be regarded as a flip;
see Koll\'ar et al.\ \cite[2.26]{Flips} on
the minimal model program for non-$\Q$-factorial varieties. For $b>1$, we are
in the more familiar case where $X_0$ is $\Q$-factorial and the flipping
contraction $X_0\arrow Z_0$ contracts a subvariety of codimension
greater than 1 in $X_0$.)

Thus, for $a>b$, we have a flip $X_0\dash Y_0$ which deforms to an isomorphism
$X_t\cong Y_t$ for $t\neq 0$.
This implies that the nef cone of $X_0$ is not the whole
nef cone of $X_t$, as we want. Indeed, an ample divisor on $Y$
restricts to an ample divisor on $Y_0$ and also on $Y_t\cong X_t$
for $t\neq 0$, but the corresponding divisor on $X_0$ is not ample
because $X_0\dash Y_0$ is not an isomorphism. For $a=b$, the flop
$X_0\dash Y_0$ is the Mukai flop, which replaces $\P^a\subset X_0$
with normal bundle $T^*\P^a$ by another copy of $\P^a$ \cite{Mukai}.
The Mukai flop is well known to have
the property here: it deforms into an isomorphism
\cite[Theorem 4.6]{Huybrechts}. The flip $X_0\dash Y_0$ for $a>b$
was considered by Kawamata \cite[Example 4.2]{Kawamataext}.

For clarity, we first describe the example in a neighborhood
of the subvariety being flipped. The flip $X\dash Y$ will be the simplest flip
between smooth varieties, with $X$ the total space of the
vector bundle $O(-1)^{\oplus b+1}$ over $\P^a$ and $Y$
the total space of the vector bundle $O(-1)^{\oplus a+1}$
over $\P^b$. These spaces can also be described in terms of linear
algebra. Let $V=A^{b+1}$ and $W=A^{a-b}$, viewed as vector spaces,
so that $V\oplus W$ has dimension $a+1$.
In this local picture, the image of the flipping contraction $X\arrow Z$
is the affine variety
$$Z=\{f=(f_1,f_2):V\arrow V\oplus W \text{ linear, } \rank(f)\leq 1\}.$$
And $X$ and $Y$ are the two obvious resolutions of $Z$, depending
on whether we look at the image or the kernel of the linear map $f$:
$$X=\{(f,L): L\subset V\oplus W \text{ a line, } f:V\arrow L
\text{ linear}\}$$
and
$$Y=\{(f,S): S\subset V \text{ of codimension 1, } f:V/S\arrow V\oplus W
\text{ linear}\}.$$

We map $Z$ to $A^1$ by the trace $t=\tr(f_1)$; composing gives
morphisms $X\arrow A^1$ and $Y\arrow A^1$. (For this local picture,
$t$ could be replaced by any sufficiently general
function on $Z$ vanishing at the origin.)
The fiber
$Y_0$ over the origin is smooth, since we can describe it
as $\{(f,S):S\subset V \text{ of codimension 1, } f:V/S\arrow S\oplus W
\text{ linear}\}$, which exhibits $Y_0$ as the total space of the vector
bundle $T^*\P^b\oplus O(-1)^{\oplus a-b}$ over $\P^b$. By contrast,
$X_0$ has singular set of codimension $2b+1$ (here $X_0$ has dimension
$a+b$, so $X_0$ is actually smooth in the flop case, $a=b$).
Explicitly,
if $x_0,\ldots,x_a$ are homogeneous coordinates
on $\P^a$ and $\xi_0,\ldots,\xi_b$
are fiber coordinates on $O(-1)$, then the morphism $X\arrow A^1$
is given by $t=\sum_{i=0}^b x_i\xi_i$. We see that the singularity
of $X_0$ is locally a node on a $(2b+1)$-fold times a smooth variety.
In particular, $X_0$ is terminal and Gorenstein for any $b\geq 1$,
and it is $\Q$-factorial for $b\geq 2$, as we want.

Now let us compactify $X$, $Y$ and $Z$ over $A^1$. We will use
the same names for the compactified varieties. Namely, define
$$Z=\{([f_1,f_2,g],t): (f_1,f_2):V\arrow V\oplus W \text{ of rank}
\leq 1, g\in A^1, t\in A^1, \tr(f_1)=tg\}.$$
(The brackets $[f_1,f_2,g]$ mean that $f_1,f_2,g$ are not all zero
and we only consider them up to a common scalar factor.) Taking $g=1$
gives the open subset of $Z$ we considered before. The number
$t$ gives the projective morphism $Z\arrow A^1$ we want.
The varieties
$X\arrow Z$ and $Y\arrow Z$ are defined by
\begin{multline*} X=\{([f_1,f_2,g],L,t): L\subset V\oplus W \text{ a line,}\\ 
(f_1,f_2):V\arrow L\subset V\oplus W, g\in A^1, t\in A^1, \tr(f_1)=tg\}
\end{multline*}
and
\begin{multline*}
Y=\{([f_1,f_2,g],S,t): S\subset V \text{ of codimension 1,} \\
(f_1,f_2):V/S\arrow V\oplus W, g\in A^1, t\in A^1, \tr(f_1)=tg\}.
\end{multline*}
Write $X_0,Y_0,Z_0$ for the fibers of $X,Y,Z$ over 0 in $A^1$.
The varieties $X$ and $Y$ are isomorphic
outside $\P^a\subset X_0$ and $\P^b\subset Y_0$, which are contracted
to a point in $Z_0$. In fact, the projective morphisms $X\arrow A^1$
and $Y\arrow A^1$ are products over $A^1-0$,
with all fibers isomorphic to $\P^a\times \P^b$, by mapping
to $(L,S)\in \P^a\times \P^b$. Here $X$ and $Y$  are smooth varieties,
but $X_0$ is singular
in codimension $2b+1$ (with its singular set contained in
$\P^a\subset X_0$). By the local calculation above,
$X_0$ is terminal and Gorenstein for any $b\geq 1$,
and it is $\Q$-factorial for $b\geq 2$. The flipped variety $Y_0$
is smooth; in fact, it is the projective bundle $P(TP^b\oplus
O(1)^{\oplus a-b}\oplus O)$ over $\P^b$ (the space of hyperplanes in the given
vector bundle). In particular, $Y_0$ has Picard number 2.

Finally, let us show that $X_0$ is Fano for $a>b$. (For $a=b$, $X_0$
and $Y_0$ are weak Fano, since we compute that the line bundles
$-K_{X_0}$ and $-K_{Y_0}$ are semi-ample and give the birational contractions
$X_0\arrow Z_0$ and $Y_0\arrow Z_0$.) The Picard number of $X_0$
is 2 (when $b>1$, this is immediate from the small $\Q$-factorial
modification $X_0\dash Y_0$). Therefore the nef cone
of $X_0$ has exactly two extremal rays, which we can
compute by exhibiting two nontrivial contractions of $X_0$.
By construction, we have  $X_0\arrow \P^a$ and $X_0\arrow Z_0$.
Both contractions are $K$-negative; in the first case, this uses
that $X_0\arrow \P^a$ is a projective bundle over
the generic point of $\P^a$. Therefore $X_0$ is Fano, as we want. \qed

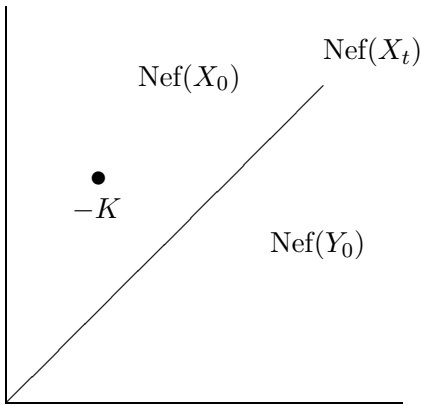
\begin{picture}(200,170)
\put(0,0){\line(1,0){150}}
\put(0,0){\line(0,1){150}}
\put(0,0){\line(1,1){120}}
\put(100,57){$\Nef(Y_0)$}
\put(50,120){$\Nef(X_0)$}
\put(120,130){$\Nef(X_t)$}
\put(25,70){$-K$}
\put(35,85){\circle*{5}}
\end{picture}

\section{Jumping of the nef cone in dimension 4}

\begin{theorem}
\label{four}
There is a variety $X$
and a flat projective morphism $t:X\arrow A^1$ with the following
properties. The fibers $X_t$ for $t\neq 0$ are isomorphic
to the blow-up of $\P^4$ along a line,
and the nef cone of $X_0$ is a proper subset 
of the nef cone of $X_t$.
The fiber $X_0$ is a $\Q$-factorial terminal Fano 4-fold.
\end{theorem}

Again, this shows the optimality of de Fernex-Hacon's results,
as discussed in the introduction.
Note that $X_0$ cannot be Gorenstein,
and so these examples are unavoidably more complicated than those in
Theorem \ref{jump}. Nonetheless, the 5-fold flip $X\dash Y$
that comes up is closely analogous to the simplest 3-fold flip
\cite[Example 2.7]{KMbook}.

{\bf Proof. }The plan of the example is the same as for
Theorem \ref{jump}. We start with a 5-fold flip $X\dash Y$
that transforms a weighted projective plane $P(1,1,2)$
into a $\P^2$, with flipping
contraction $X\arrow Z$. The flip we use can be defined as a quotient
of the 5-fold flop $(\P^2,O(-1)^{\oplus 3})\dash (\P^2,O(-1)^{\oplus 3})$
by the group $\Z/2$ (as in the simplest
example of a 3-fold flip \cite[Example 2.7]{KMbook}).
In this case (unlike in Theorem \ref{jump}),
$X$ will be singular; it has one singular point in the subvariety $P(1,1,2)$,
and $X$ is not Gorenstein at that point. We consider a general
function $t:Z\arrow A^1$ vanishing at the singular point of $Z$.
We have to check that the 4-fold $X_0:=t^{-1}(0)\subset X$
is $\Q$-factorial and terminal. With that done, we have a ``Mukai-type''
flip of 4-folds $X_0\dash Y_0$, meaning one that deforms into
an isomorphism $X_t\cong Y_t$. Finally, we compactify
$X$ over $A^1$, and check that we can arrange
for $X_0$ to be a Fano variety.

For brevity, we will define a compactified family $X\arrow A^1$
without describing the local picture in detail. First define
$$W=\{ ([x_1,x_2,x_3,y_2,y_3,g],t)\in P(1^5,2)\times A^1:
x_1^2+x_2y_2+x_3y_3=tg\}.$$
(A reference on weighted projective spaces is Dolgachev \cite{Dolgachev}.)
Then $t:W\arrow A^1$ is a flat projective morphism with fibers isomorphic
to $\P^4$ outside $0$ in $A^1$ (since we can solve the equation for $g$
when $t$ is not zero). It is a smooth morphism outside
the one singular point of $W$, $([0,0,0,0,0,1],0)$, where
$W$ has a quotient singularity $A^5/\{ \pm 1\}$. The fiber
$W_0$ over 0 in $A^1$
is the quotient by $\Z/2$ of a quadric 4-fold in $\P^5$ with a node,
and so $W_0$ has Picard number 1.

The singularity
of $W_0$ (the quotient of a 4-fold node by $\pm 1$)
is obtained by contracting a smooth quadric 3-fold $D$
with normal bundle $O(-2)$ in a smooth 4-fold $M$. Since $D$ has Picard
number 1, $W_0$ is $\Q$-factorial. We compute that $\pi:M\arrow W_0$ has
$K_{M}=\pi^*K_{W_0}+(1/2)D$.
Since the discrepancy $1/2$ is positive, $W_0$ is terminal.

Let $X$ be the blow-up of $W$ along the subvariety
$\P^1\times A^1=\{ ([0,0,0,y_2,y_3,0],t) \}$ of $W$. Since
this is contained in the smooth locus of $W\arrow A^1$, the morphism
$X\arrow A^1$ is smooth outside one point. 
Clearly the fibers of $X\arrow A^1$ outside $0$ in $A^1$ are isomorphic
to the blow-up of $\P^4$ along a line. 
Since the fiber $X_0$ is obtained by blowing up $W_0$ at a $\P^1$
in the smooth locus of $W_0$,
$X_0$ is $\Q$-factorial and terminal, and it has Picard number 2.

We now define a flipping
contraction $X\arrow Z$ over $A^1$ which contracts a surface $P(1,1,2)$
in the special fiber $X_0$ to a point. Let
$$Z=\{ ([u_{11},u_{12},u_{13},u_{21},\ldots,u_{33},g],t)\in
P(1^3,2^7)\times A^1: \rank (u_{ij})\leq 1,
u_{11}^2+u_{22}+u_{33}=tg\}.$$
Here we view $(u_{ij})$ as a $3\times 3$ matrix to define its rank.
We define the contraction $X\arrow Z$ as the following
rational map $W\dash Z$, which becomes a morphism after blowing
up the subvariety $\P^1\times A^1$ of $W$:
$$([x_1,x_2,x_3,y_2,y_3,g],t)\mapsto
  ([x_1,x_2,x_3,x_1y_2,x_2y_2,x_3y_2,x_1y_3,x_2y_3,x_3y_3,g],t).$$
The resulting contraction $X\arrow Z$ contracts only a surface
$S\cong P(1,1,2)$, the birational transform \cite[Notation 0.4]{KMbook}
of the surface $\{ ([0,0,0,y_2,y_3,g],0)\}$ in $W$.

It is straightforward to check that the contraction $X\arrow Z$
is $K_X$-negative; locally, it is a terminal toric flipping contraction,
corresponding to a relation
$e_1+e_2+2e_3=e_4+e_5+e_6$ in $N\cong \Z^5$ in Reid's notation
\cite[Proposition 4.3]{Reid}.
The toric picture shows that the flip $Y\arrow Z$ of $X\arrow Z$
replaces the surface $P(1,1,2)$ in the 5-fold $X$ with a copy of $\P^2$.

Since $X_0$ has Picard number 2 and has two nontrivial $K$-negative
contractions (to $W_0$ and $Z_0$), $X_0$ is Fano.
Because $X\arrow Z$ contracts
only a surface in $X_0$, $X$ and $Y$ are
the same except over 0 in $A^1$.
That completes the proof: $X_0$ is a $\Q$-factorial terminal
Fano 4-fold with a flip $X_0\dash Y_0$ that deforms into
an isomorphism $X_t\cong Y_t$. It follows that the nef cone
of $X_0$ is a proper subset of that of $X_t$, because an ample
divisor on $Y$ is ample on $X_t\cong Y_t$ but not on $X_0$. \qed

\section{The divisor class group of a Fano variety}

\begin{theorem}
\label{top}
Let $X$ be a complex projective $n$-fold with rational singularities
such that $H^1(X,O)=H^2(X,O)=0$.
Then $\Pic(X)\arrow H^2(X,\Z)$ and $\Cl(X)\arrow
H_{2n-2}(X,\Z)$ are isomorphisms.
\end{theorem}

The theorem applies, for example,
if there is an effective $\R$-divisor $\Delta$
with $(X,\Delta)$ a dlt Fano pair. 
The statement on the Picard group
is well known, but it is harder for the divisor
class group, even for terminal Fano varieties.
When $X$ has dimension 3 and the singularities are isolated,
the statement on the divisor class group
was proved by Namikawa-Steenbrink \cite[Theorem 3.2]{NS}.

{\bf Proof. }To say that $(X,\Delta)$ is a dlt Fano pair means
that the pair $(X,\Delta)$ is dlt \cite[Definition 2.37]{KMbook}
and that $-(K_X+\Delta)$ is ample. In that case, 
$X$ has rational singularities \cite[Theorem 5.22]{KMbook},
and $H^i(X,O)=0$ for $i>0$ by Kodaira vanishing
\cite[Theorem 2.42]{Fujinobook}.

Let $X$ be a projective variety with rational singularities
such that $H^1(X,O)=H^2(X,O)=0$. By the exponential
sequence, $\Pic(X)\arrow H^2(X,\Z)$ is an isomorphism.
Let $\pi:\tx\arrow X$ be a resolution of singularities, which
we can assume to be an isomorphism over the smooth locus $X^{\sm}$ and with
the inverse image of the singular set a divisor $E=\cup E_i$ having simple
normal crossings. By definition of rational singularities,
we have $\pi_*O_{\tx}=O_X$ and $R^j\pi_*O_{\tx}=0$ for $j>0$.
We deduce that
$H^1(\tx,O)=H^2(\tx,O)=0$, and so $\Pic(\tx)$ maps isomorphically
to $H^2(\tx,\Z)$. Equivalently, $\Cl(\tx)$ maps isomorphically 
to $H_{2n-2}(\tx,\Z)$.

We have localization sequences in Chow groups and Borel-Moore homology,
using that $\tx-E\cong X^{\sm}$:
$$\begin{CD}
\oplus \Z E_i @>>> \Cl(\tx)@>>> \Cl(X^{\sm}) @>>> 0\\
@VV{\cong}V @VV{\cong}V @VVV @VVV \\
\oplus \Z E_i @>>> H_{2n-2}(\tx,\Z) @>>> H_{2n-2}^{BM}
(X^{\sm},\Z) @>>> H_{2n-3}(E,\Z).
\end{CD}$$
The diagram implies that $\Cl(X^{\sm})\arrow H_{2n-2}^{BM}(X^{\sm},\Z)$
is injective. Equivalently, $\Cl(X)\arrow H_{2n-2}(X,\Z)$ is injective.

Namikawa showed that for any variety $X$ with rational
singularities, all 2-forms on the smooth locus extend to 2-forms
on a resolution $\tx\arrow X$ \cite[Theorem A.1]{Namikawa}.
In our situation,
we deduce that every 2-form on $\tx$ with log poles on the divisor
$E$ extends to a 2-form on $\tx$, which must be zero by the Hodge symmetry
$h^{pq}=h^{qp}$ from $H^2(\tx,O)=0$.

By definition of Deligne's Hodge filtration on $H^2(X^{\sm},\C)$, the
vanishing of 2-forms on $\tx$ with log poles on $E$
means precisely that $F^2H^2(X^{\sm},\C)=0$ \cite[Theorem 3.2.5]{Deligne}.
The mixed Hodge structure on $H^2$ of any smooth variety
has Hodge numbers $h^{pq}$ equal to zero 
unless $p\leq 2$, $q\leq 2$, and $p+q\geq 2$
\cite[Corollary 3.2.15]{Deligne}. We have shown that $h^{20}=h^{21}=h^{22}=0$;
using Hodge symmetry, it follows that
the only nonzero Hodge number for $H^2(X^{\sm},\C)$ is $h^{11}$. In particular,
$H^2(X^{\sm},\C)$ is all in weight 2, which means that $H^2(\tx,\Q)
\arrow H^2(X^{\sm},\Q)$ is surjective \cite[Corollary 3.2.17]{Deligne}.
Here $\Cl(\tx)$ is isomorphic to $H^2(\tx,\Z)$ and therefore maps onto
$H^2(X^{\sm},\Z)=H_{2n-2}(X,\Z)$
after tensoring with the rationals. It follows that
$\Cl(X)$ maps onto $H_{2n-2}(X,\Z)$ after tensoring
with the rationals.

To show that $\Cl(X)$ maps onto $H_{2n-2}(X,\Z)$, we observe
that the cokernel of $\Cl(X)\arrow H_{2n-2}(X,\Z)$ is torsion-free
for every normal algebraic variety $X$. Indeed, it is equivalent
to show that $\Pic(X^{\sm})\arrow H^2(X^{\sm},\Z)$ has torsion-free cokernel.
That follows from the commutative diagram, for any positive integer $m$,
where the top row is the Kummer sequence for the etale topology:
$$\begin{CD}
H^1(X^{\sm},\Z/m)@>>> \Pic(X^{\sm}) @>>{\cdot m}> \Pic(X^{\sm})
@>>> H^2(X^{\sm},\Z/m)\\
@VV{\cong}V @VVV @VVV @VV{\cong }V \\
H^1(X^{\sm},\Z/m)@>>> H^2(X^{\sm},\Z) @>>{\cdot m}> H^2(X^{\sm},\Z)
@>>> H^2(X^{\sm},\Z/m).\\
\end{CD}$$
We have shown that $\Cl(X)\arrow H_{2n-2}(X,\Z)$ is an isomorphism. \qed

\section{Deforming Weil divisors}

\begin{theorem}
\label{class}
Let $X_0$ be a complex projective variety with rational singularities
such that $H^1(X_0,O)=H^2(X_0,O)=0$.
Suppose that $X_0$ is
smooth in codimension 2 and $\Q$-factorial in codimension 3.
Then the divisor class group is unchanged under nearby deformations
of $X_0$. More precisely, given a deformation $X\arrow (T,0)$
of $X_0$ over a smooth curve $T$,
there is an etale morphism of curves
$(T',0)\arrow (T,0)$ such that the class group of $X_{T'}$ maps 
split surjectively to
the class group of $\Cl(X_t)$ for all $t\in T'$, and all these
surjections have the same kernel.
\end{theorem}

The assumption that $X_0$ has rational singularities and $H^1(X_0,O)
=H^2(X_0,0)=0$
holds, for example,
if there is an effective $\R$-divisor $\Delta$
with $(X_0,\Delta)$ a dlt Fano pair. 

In Theorem \ref{class},
the divisor class group means the group of Weil divisors
(with integer coefficients) modulo linear equivalence. For
$\Q$-factorial terminal Fanos, deformation invariance
of the class group was proved by
de Fernex and Hacon \cite[Proposition 6.4, Lemma 7.2]{DH}, building on
the work of Koll\'ar and Mori \cite[Proposition 12.2.5]{KMflip}.

The assumptions on singularities
are explained by some simple examples. A quadric
surface $X_0$ with a node is $\Q$-factorial but singular in codimension 2,
and it deforms to a smooth quadric surface $X_t$; then $\Cl(X_0)\cong \Z$
and $\Cl(X_t)\cong\Z^2$. A quadric 3-fold $X_0$ with a node is
smooth in codimension 2 but not $\Q$-factorial in codimension 3,
and it deforms to a smooth quadric 3-fold $X_t$; then $\Cl(X_0)\cong \Z^2$
and $\Cl(X_t)\cong \Z$. Generalizing the last case, Sturmfels,
Batyrev, and others found many 
non-toric Fanos such as Grassmannians which degenerate to
terminal Gorenstein toric Fano varieties $X_0$ \cite{Batyrev}.
In their examples,
$\Pic(X_0)\cong \Pic(X_t)\cong \Z$, but $X_0$ is not $\Q$-factorial
because of nodes in codimension 3, whereas $X_t$ is smooth. So
$\Cl(X_0)$ is bigger than $\Cl(X_t)\cong \Z$. These examples suggest
the problem of finding  ``minimal'' assumptions which imply
that the Picard group of a Fano variety, rather than the divisor
class group, is constant under deformation.

Finally, Theorem \ref{class} does not hold
if the assumption that $X_0$ has rational singularities
is weakened to $X_0$ Cohen-Macaulay. An example is provided
by taking $X_0$ to be the projective cone over a smooth quartic
surface $Y_0$ with Picard number 1, and deforming it to projective cones $X_t$
over other smooth quartic surfaces $Y_t$.
Then $X_0$ is Cohen-Macaulay, smooth in codimension 2, $\Q$-factorial,
and $H^1(X_0,O)=H^2(X_0,O)=0$. (It is also a log-canonical Fano 3-fold.)
But $\Cl(X_t)\cong \Pic(Y_t)$ can be non-constant even on an analytic
neighborhood of $0\in T$. On the good side, if the assumption
that $X_0$ has rational singularities in Theorem \ref{class}
is weakened to $X_0$ Cohen-Macaulay, then the proof at least shows
that there is an etale morphism of curves $(T',0)\arrow (T,0)$
such that the class group of $X_{T'}$ maps split surjectively
to the class group of $X_0$.

{\bf Proof. }We know that $\Cl(X_0)$ is finitely generated by
the more precise Theorem \ref{top}, $\Cl(X_0)=H_{2n-2}(X_0,\Z)$.
Let us prove that $\Cl(X_R)\arrow \Cl(X_0)$ is injective, where
$R$ is the henselization of the local ring of $T$ at 0
\cite[section I.4]{Milne} and $X_R$ denotes the pullback
of $X\arrow T$ via $\Spec(R)\arrow T$. 
(We will briefly replace $T$ by a ramified covering later in this proof,
but the same argument will apply.) Namely, by Griffith's local
Lefschetz theorem for divisor class groups, since $X_0$ is
Cohen-Macaulay and smooth in codimension 2, the restriction map
$\Cl(O_{X_R,x})\arrow \Cl(O_{X_0,x})$ is injective for every point
$x\in X_0$ \cite[Proposition 3.5]{Griffith}. Therefore a
Weil divisor class $D$ on $X_R$ which restricts to 0 in $\Cl(X_0)$ is
trivial near each point of $X_0$. This means that $D$ is a Cartier divisor,
and its restriction to $X_0$ is linearly equivalent to zero. But
$\Pic(X_R)\arrow \Pic(X_0)$ is injective because $H^1(X_0,O)=0$
\cite[Corollary 4.6.4]{EGA3}.
So $D$ is linearly equivalent to zero on $X_R$. 
We have shown that $\Cl(X_R)\arrow \Cl(X_0)$ is injective.

Let $K$ be the quotient field of the henselian local ring $R$.
Geometrically,
$X_K$ is obtained from $X_R$ by removing the special fiber $X_0$, which
is linearly equivalent to zero as a divisor in $X_R$,
and so the restriction map
$\Cl(X_R)\arrow \Cl(X_K)$ is an isomorphism. The algebraic closure
$\Kb$ of $K$ has Galois group $\Gal(\Kb/K)\cong \Zh$ (the completed
fundamental group of a ball in $T$ minus the point 0).
I claim that $\Cl(X_R)=\Cl(X_K)$ maps isomorphically to the Galois
invariants $\Cl(X_{\Kb})^{\Zh}$. To prove this,
let $U$ be the smooth locus of $X_K$;
since $X_K$ is normal, it suffices to show that $\Pic(U)\arrow \Pic(U_{\Kb})^
{\Zh}$ is an isomorphism. Consider the Hochschild-Serre spectral
sequence for etale cohomology,
$$H^p(K,H^q_{\et}(U_{\Kb},O^*))\imp H^{p+q}_{\et}(U,O^*).$$
Here $H^0_{\et}(U_{\Kb},O^*)\cong \Kb^*$ since $X_K$ is normal, connected,
and proper over $K$.
So $H^1(K,H^0_{\et}(U_{\Kb},O^*))=H^1(K,\Kb^*)=0$
by Hilbert's theorem 90. Also,
$H^2(K,H^0_{\et}(U_{\Kb},O^*))$ is equal to
the Brauer group $H^2(K,\Kb^*)$, which is zero because $K$ is the quotient
field of a henselian discrete valuation ring with algebraically
closed residue field \cite[Example III.2.2]{Milne}.
Therefore the spectral sequence gives that
$\Pic(U)=H^1_{\et}(U,O^*)$ maps isomorphically to $\Pic(U_{\Kb})^{\Zh}
=H^0(K,H^1_{\et}(U_{\Kb},O^*))$. Thus we have shown that $\Cl(X_R)=\Cl(X_K)$
maps isomorphically to $\Cl(X_{\Kb})^{\Zh}$.

Notice that the class group
$\Cl(X_{\Kb})$ is isomorphic to the class group of $X_t$
for very general points $t\in T$ (that is, for all but countably
many points in the curve $T$). Indeed, the morphism $X\arrow T$
is defined over some finitely generated subfield $F$ of $\C$,
$X_F\arrow T_F$. For very general points $t\in T(\C)$, the morphism
$t:\Spec(\C)\arrow T\arrow T_F$ factors through a geometric
generic point, $\overline{F(T)}$ (that is, it is transcendental over $F$).
So both $X_t$ and $X_{\Kb}$ are extensions
of $X_{\overline{F(T)}}$ from one algebraically closed field to another.
Since the class group of $X_t$ is finitely generated, this implies
that the class group of $X_t$ (for $t$ very general) is isomorphic
to that of $X_{\Kb}$.

I claim that the Galois group $\Zh$ acts trivially on $\Cl(X_{\Kb})$.
(In geometric terms, this means that the monodromy $\Z$ around 0 in $T$
acts trivially on $\Cl(X_t)=H_{2n-2}(X_t,\Z)$ for $t\neq 0$ near 0.)
Our earlier proof that $\Cl(X_R)=\Cl(X_K)$
maps injectively to $\Cl(X_0)$ applies to any finite extension field
$E$ in place of $K$, to show that
the specialization map $\Cl(X_E)\arrow \Cl(X_0)$ is injective.
By taking a direct limit, it follows that
the specialization map  $\Cl(X_{\Kb})\arrow \Cl(X_0)$ is injective.
But this specialization
map factors through the coinvariants of the Galois group
$\Gal(\Kb/K)=\Zh$ on $\Cl(X_{\Kb})$. Therefore, $\Gal(\Kb/K)=\Zh$
acts trivially on $\Cl(X_{\Kb})$, as we wanted.
Together with the argument two paragraphs back, this shows that
$\Cl(X_R)=\Cl(X_K)$ maps isomorphically to $\Cl(X_{\Kb})$.

Since $R$ is the henselian local ring of $T$ at 0, this implies that 
after replacing $T$ by
some algebraic curve $T'$ with an etale morphism $(T',0)\arrow
(T,0)$, the class group $\Cl(X)$ maps onto $\Cl(X_t)$ for
very general $t\in T$. We can assume (after replacing $T$ by
a Zariski open neighborhood of 0)
that $X\arrow T$ is a topological fibration over $T-0$. 
Since we have a topological interpretation of the divisor
class group as $\Cl(X_t)=H_{2n-2}(X_t,\Z)$, the class group
$\Cl(X)$ maps onto $\Cl(X_t)$ for all $t\neq 0$ in $T$.

It remains to show that $\Cl(X_R)$ maps onto $\Cl(X_0)$.
Let $D$ be a Weil divisor on $X_0$. The sheaf
$O(D)$ restricts to a line bundle $L$ on the smooth locus $U_0$ of $X_0$.
Write $jU_0$ for the subscheme of $X$ defined locally by the $j$th
power of a function defining $U_0$. We have an exact sequence
$$H^1(U_0,O)\arrow \Pic((j+1)U_0)\arrow \Pic(jU_0)\arrow H^2(U_0,O),$$
because the normal bundle of $X_0$ in $X$ is trivial.
Because $X_0-U_0$ has codimension at least 3 and $X_0$ is Cohen-Macaulay,
the restriction map $H^j(X_0,O)\arrow H^j(U_0,O)$ is bijective
for $i<2$ and injective for $i=2$ \cite[Proposition III.3.3]{SGA2}. Since
$H^1(X_0,O)=0$, that implies that $H^1(U_0,O)=0$. Moreover,
$U_0$ is contained in some $\Q$-factorial
open subset $V_0$ whose complement has codimension at least 4 in $X_0$.
By the same arguments, $H^1(V_0,O)=H^2(V_0,O)=0$. 

Because $V_0$ is $\Q$-factorial, there is a positive integer $m$
such that the line bundle $L^{\otimes m}$ on $U_0$ extends
to a line bundle $M$ on $V_0$.
By the exact sequence
we wrote for $U_0$, $M$ extends uniquely to a line
bundle $M$ on the formal scheme $V_0^{\wedge}$ (the completion of $V_0$ in $X$).
It follows that the obstruction in $H^2(U_0,O)$ to extending $L$
from $U_0$ to $2U_0$ is killed by $m$ and hence is zero, because
$H^2(U_0,O)$ is a complex vector space. Because $H^1(U_0,O)=0$,
the extension of $L$ to $2U_0$ is unique, and repeating
the argument shows that $L$ extends to a line bundle
on the formal scheme $U_0^{\wedge}$, the completion of $U_0$ in $X$.

Let $j:U_0^{\wedge}\arrow X_0^{\wedge}$ be the inclusion of formal schemes,
where $X_0^{\wedge}$ is the completion of $X_0$ in $X$. 
Because $X_0$ is Cohen-Macaulay and the codimension of $U_0$ in $X_0$
has codimension at least 3, Grothendieck's existence theorem
shows that $E:=j_*L$ is a coherent sheaf on the formal scheme $X_0^{\wedge}$
\cite[Theorem IX.2.2]{SGA2}. Because $X\arrow T$ is a projective
morphism, $E$ is algebraizable \cite[Corollary III.5.1.6]{EGA3};
that is, $E$ can be viewed as a coherent sheaf
on the scheme $X_{R^{\wedge}}$, where
$R^{\wedge}$ is the completed local ring of the curve $T$ at the point $0$.
Because the restriction of $E$ to $U_0\subset X_0$ is a line bundle,
$E$ is a line bundle on $X_{R^{\wedge}}-S$ for some closed subset
$S\subset X_{R^{\wedge}}$ of codimension at least 3. So, after
replacing the Weil divisor
$D$ on $X_0$ we started with by a linearly equivalent divisor,
$D$ extends to a Weil divisor $D$ on $X_{R^{\wedge}}$.

The relative Hilbert scheme of codimension-1 subschemes of fibers
of $X\arrow T$ with given Hilbert polynomial $p$, $\Hilb^{p}(X/T)$,
is projective over $T$. Write $D=D_1-D_2$ where $D_1$ and $D_2$
are effective divisors on $X_{R^{\wedge}}$. A projective morphism which has a section
over the completed local ring of $T$ at 0 also has a section
over the henselian local ring $R$ of $T$ at 0. Therefore, $D$ extends
to a Weil divisor on $X_R$. We have shown that $\Cl(X_R)
\arrow \Cl(X_0)$ is surjective.

We showed at the
start of the proof that
$\Cl(X_R)\arrow \Cl(X_0)$ is injective, and so
this restriction map is an isomorphism. (We also know
that $\Cl(X_R)\arrow \Cl(X_{\Kb})$ is an isomorphism.) But $\Cl(X_R)$
is the direct limit of the class groups of all pullbacks of
$X$ by etale morphisms $(T',0)\arrow (T,0)$. Therefore, after
some such pullback, the homomorphism $\Cl(X)\arrow \Cl(X_0)$ becomes
split surjective. By what we have shown, the kernels of the restriction
maps $\Cl(X)\arrow \Cl(X_0)$ and $\Cl(X)\arrow \Cl(X_t)$ for $t\neq 0$
in $T$ are both equal to the kernel of $\Cl(X)\arrow \Cl(X_R)$.
\qed

\section{Deformations of toric Fano varieties}

\begin{theorem}
\label{toric}
A toric Fano variety which is smooth in codimension 2
and $\Q$-factorial in codimension 3 is rigid.
\end{theorem}

This strengthens Bien-Brion's theorem that smooth toric Fano
varieties are rigid \cite{BB} and de Fernex-Hacon's theorem
that $\Q$-factorial terminal toric Fano varieties are rigid
\cite{DH}. (Note that a toric variety is $\Q$-factorial if and only
if the corresponding fan is simplicial.) The assumptions
on singularities cannot be omitted. First,
a quadric surface with a node is a
toric Fano variety that is not rigid. Likewise, many toric
Fano varieties with codimension-3 node singularities are not rigid,
by the examples after Theorem \ref{class}.
The varieties in Theorem \ref{toric} are klt (like every 
$\Q$-Gorenstein toric variety) but need not be terminal or canonical.

{\bf Proof. }For any normal variety $X$ and $j\geq 0$, the double dual
$(\Omega^j_X)^{**}$
is the reflexive sheaf that extends the vector bundle $\Omega^j_U$
on the smooth locus $U$ of $X$; that is, $(\Omega^j_X)^{**}=
f_*\Omega^j_U$ where
$f:U\arrow X$ is the inclusion.
Danilov's vanishing theorem says that for any ample line
bundle $O(D)$
on a projective toric variety $X$, $H^i(X,(\Omega^j\otimes O(D))^{**})=0$
for all $i>0$ and $j\geq 0$ \cite[Theorem 7.5.2]{Danilov}. More generally,
for any ample Weil divisor $D$ on a projective toric variety $X$,
Musta\textcommabelow{t}\u{a}
showed that $H^i(X,(\Omega^j\otimes O(D))^{**})=0$
for all $i>0$ and $j\geq 0$
\cite[Proposition 2.3]{Mustata}.
Fujino gave a simple proof of this vanishing using the action
of the multiplicative monoid
of natural numbers on a toric variety \cite[Proposition 3.2]{Fujinomult}.

Define the tangent sheaf of a normal variety $X$ of dimension $n$
as $\T^0_X=\Hom(\Omega^1_X, O_X)$. This is the reflexive sheaf
that agrees with $\Omega^{n-1}_X\otimes (K_X)^*$ on the smooth locus
of $X$, and so we have $\T^0_X=(\Omega^{n-1}\otimes (K_X)^*)^{**}$.
By the Danilov-Musta\textcommabelow{t}\u{a}
vanishing theorem, for any toric Fano variety
$X$, we have $H^i(X,\T^0_X)=0$ for all $i>0$. Since
$H^1(X,\T^0_X)=0$, all locally trivial deformations
of $X$ are trivial \cite[Proposition 1.2.9]{Sernesi}.
(Roughly speaking, a deformation
is locally trivial if it does not change the singularities of $X$.)

Finally, Altmann showed that any $\Q$-Gorenstein affine toric variety
which is smooth in codimension 2 and $\Q$-factorial in codimension 3
is rigid \cite[Corollary 6.5(1)]{Altmann}.
Thus, for $X$ a toric Fano variety which
is smooth in codimension 2 and $\Q$-factorial in codimension 3,
all deformations are locally trivial. By the previous paragraph,
all deformations of $X$ are in fact trivial. \qed

\section{Flatness of Cox rings}

\begin{theorem}
\label{flat}
Let $X_0$ be a terminal Fano variety over the complex numbers
which is $\Q$-factorial
in codimension 3. Then, for any deformation $X\arrow (T,0)$ of $X_0$
over a smooth curve,
the Cox rings $\Cox(X_t)$ form a flat family over an etale neighborhood
of 0 in $T$. That is, after replacing $X$ by its pullback via
some etale morphism of curves $(T',0)\arrow (T,0)$, every divisor class
on each $X_t$ extends to $X$, and
$h^0(X_t,O(D))$ is independent of $t$ for every Weil divisor $D$ on $X$.
\end{theorem}

Here by the Cox ring of $X_t$
we mean the sum over all classes of Weil divisors, $\oplus_{D\in \Cl(X_t)}
H^0(X_t,O(D))$. For $X_0$ $\Q$-factorial rather than $\Q$-factorial
in codimension 3, the theorem was proved by de Fernex and Hacon
\cite[Theorem 1.1, Proposition 6.4]{DH}. The proof of Theorem \ref{flat}
combines their results with Theorem \ref{class} on deformation invariance
of the divisor class group.

The statement about flatness in Theorem \ref{flat} can be explained
in more detail. Let $X\arrow (T,0)$
be a deformation over an affine curve,
as in the theorem. Theorem \ref{class}
shows that after replacing $X$
by its pullback by some etale morphism of curves
$(T',0)\arrow (T,0)$, we can assume that $\Cl(X)\arrow \Cl(X_t)$ is split
surjective for all $t$ in the new curve $T$, and all these surjections
have the same kernel.
Let $A\subset \Cl(X)$
be a subgroup which maps isomorphically to $\Cl(X_t)$ for all 
$t\in T$. By choosing divisors whose classes generate the abelian group
$A$, we can define a ring structure on $\oplus_{D\in A} H^0(X,O(D))$,
which depends on choices because the group $A$ may have
torsion. 
Then Theorem \ref{flat} shows that this ring is a flat $O(T)$-algebra.
Tensoring it with the homomorphism $O(T)\arrow \C$ associated to any
point $t\in T$ gives the Cox ring of $X_t$.

{\bf Proof. }By Kawamata and Nakayama,
$X$ is terminal, after shrinking $T$ to a
Zariski open neighborhood of $0\in T$ \cite[Theorem 1.5]{Kawamataext}.
Let $\pi: W\arrow X$ be a $\Q$-factorialization,
which exists by Birkar-Cascini-Hacon-M\textsuperscript{c}Kernan 
\cite[Corollary 1.4.3]{BCHM}.
That is, $W$ is $\Q$-factorial
and terminal and $W\arrow X$ is a small projective birational morphism.
It follows that $K_W=\pi^*(K_X)$. 

We will show that $W_0\arrow X_0$ is a $\Q$-factorialization.
Since $W\arrow X$ is a small birational morphism,
$W_0\arrow X_0$ is birational. Since $X$ is terminal (in particular,
$K_X$ is $\Q$-Cartier) and $X_0$ is terminal, the pair
$(X,X_0)$ is plt \cite[Theorem 5.50]{KMbook}.
Therefore
the pair $(W,W_0)$ is plt, since $K_W=\pi^*(K_X)$ and $W_0=\pi^*(X_0)$
as a Cartier divisor. It follows that $W_0$ is normal
\cite[Proposition 5.51]{KMbook}.

Since $K_{W_0}=K_W|_{W_0}$ and $K_{X_0}=K_X|_{X_0}$, we have
$K_{W_0}=\pi^*K_{X_0}$. Since $X_0$ is terminal, it follows that
$W_0$ is terminal and the birational morphism $W_0\arrow X_0$ is small.

By Theorem \ref{class}, after replacing $T$ by a curve $T'$ with
an etale morphism $(T',0)\arrow (T,0)$,
our assumptions imply that $\Cl(X)\arrow \Cl(X_0)$ is onto.
That is, every Weil divisor on $X_0$
is the restriction of some Weil divisor $D$ on $X$. Let $D_W$
be the birational transform of $D$ on $W$. Since $W$ is $\Q$-factorial,
$D_W$ is $\Q$-Cartier, and hence so is $D_W|_{W_0}$. Clearly this
divisor pushes forward to $D|_{X_0}$, and since $W_0\arrow X_0$
is small, $D_W|_{W_0}$ must be the birational transform of $D|_{X_0}$.
Every divisor on $W_0$ is the birational transform of a divisor on $X_0$,
and so we have shown that $W_0$ is $\Q$-factorial. That is, the morphism
$W_0\arrow X_0$ is a $\Q$-factorialization. In other words, our
assumptions imply that ``$\Q$-factorialization works in families.''

Because $W_0\arrow X_0$ is small, so is $W_t\arrow X_t$ for
all $t\in T$, after shrinking $T$ around 0.
Therefore the Cox rings of these
two varieties are the same:
$\Cox(W_t)\cong \Cox(X_t)$. So, to prove flatness of the Cox rings
for the family $X\arrow T$, it suffices to prove flatness of the Cox
rings for the family $W\arrow T$. Here $W_0$ is no longer Fano, but
it is a $\Q$-factorial terminal weak Fano (meaning that $-K_{W_0}$
is nef and big), and $X_0$ is the anticanonical
model of $W_0$.

We know that the class groups of $W_t$ are independent of $t$
by Theorem \ref{class}. Our goal is to show that
for every integral Weil divisor $L$ on $X$, the restriction map
$H^0(W,O(L))\arrow H^0(W_0,O(L|_{W_0}))$ is surjective.
We can assume
that $H^0(W_0,O(L|_{W_0}))$ is not zero.

Because $\pi:W\arrow X$ is birational, there is an ample divisor $A$ on 
$X$ such that $L+\pi^*A$ is $\Q$-linearly equivalent to an effective
$\Q$-divisor $B$ on $W$. 
Choose $e>0$ small enough that
the pair $(W,eB)$ is terminal. Here
$K_{W}+eB$ is big over $X$ (just because $W\arrow X$
is birational), and so Birkar-Cascini-Hacon-M\textsuperscript{c}Kernan
showed that $(W,eB)$ has a minimal model $W'$ over $X$
\cite[Theorem 1.2]{BCHM}.
Because the morphism
$W\arrow X$ is small, so is $W'\arrow X$,
and we end up with another
small $\Q$-factorialization $\pi:W'\arrow X$ such that $K_{W'}+eB$
is nef over $X$. Equivalently, since $K_{W'}=\pi^*K_{X}$
and $B\sim_{\Q} L+\pi^*(A)$, 
$L$ becomes nef over $X$ on $W'$. Let us write $W$ instead of $W'$;
our earlier arguments still apply to this new $\Q$-factorialization
of $X$. In particular, $W_0\arrow X_0$ is a small $\Q$-factorialization,
and we have now arranged that $L$ is nef over $X$ on $W$.
Making a small $\Q$-factorial modification does not change the
space of global sections of a reflexive sheaf of rank one, and so
our goal is still to show that $H^0(W,O(L))\arrow
H^0(W_0,O(L|_{W_0}))$ is surjective.

Let $D_1$ be an effective $\Q$-divisor $\sim_{\Q} -K_W+bL$ for some
small $b>0$; this is possible because $-K_W$ is big over $T$.
Let $B$ be $1/m$
times a general divisor in $|-mK_W|$ for some $m>1$ such that $-mK_W$
is basepoint-free.
Let $D=(1-c/b)B+(c/b) D_1$
for a small $c\in (0,1)$; clearly $D\sim_{\Q} -K_W+cL$, and the pair
$(W,D)$ is klt for $c$ small enough \cite[Theorem 4.8]{Kollarsing}.
Since we arranged
that $L$ is nef on $W$ over $X$ and $-K_W=\pi^*(-K_X)$,
$D$ is nef on $W$ over $X$.

Since $-K_{X_0}$ is ample, $-K_X$ is ample over $T$ after shrinking
the curve $T$ around the point 0. So $-K_W=\pi^*(-K_X)$
is nef and big over $T$. By the cone theorem,
the cone of curves of $W$ over $T$ is rational polyhedral
\cite[Theorem 3.7]{KMbook}.  (Since $-K_W$ is nef and big over $T$,
and $W$ is $\Q$-factorial and klt, it is straightforward to construct
an effective $\Q$-divisor $E$ on $W$ such that the pair $(W,E)$ is klt
and Fano over $T$. That is the assumption we need for the cone theorem.)
Since $D$ is nef over $X$ and $-K_W$
is the pullback of $-K_X$ which is ample over $T$, $D-aK_W$
is nef over $T$ for $a$ at least some positive real number $a_0$,
since it suffices to check
that $D-aK_W$ has nonnegative intersection with each extremal
ray of the cone of curves of $W$ over $T$. Choose an $a>a_0$.

Following the idea of de Fernex-Hacon \cite[Theorem 4.5(b)]{DH},
we run a $(W,D)$-minimal model program over $T$
with scaling of $D-aK_W$:
$$W=W^0\dash W^1\dash W^2\dash\cdots .$$
(In order to run this minimal model program, we need to know
that $(K_W+D)+t(D-aK_W)$ is nef over $T$ for some $t>0$,
which is true by our choice of $a$.)
Here $D-aK_W$ is only nef over $T$, not ample over $T$ as in
\cite[Theorem 4.5(b)]{DH}, but because $W$ is weak Fano over $T$,
we still know that this minimal model program
terminates. In fact, every minimal model program on $W$
over $T$ terminates, by Birkar-Cascini-Hacon-M\textsuperscript{c}Kernan
\cite[Corollary 1.3.1]{BCHM}.

By de Fernex-Hacon \cite[Theorem 4.1]{DH}, which applies
to our minimal model program with scaling by the nef divisor $D-aK_W$,
each fiber type (resp.\ divisorial, resp.\ flipping) 
contraction of $W$ restricts
to a fiber type (resp.\ divisorial, resp.\ flipping) contraction of $W_0$.
Therefore, the given $(W,D)$-minimal model program over $T$
induces an $(W_0,D|_{W_0})$-minimal model program on $W_0$.
We assumed that the $\Q$-divisor class
$K_{W_0}+D|_{W_0}\sim_{\Q} cL|_{W_0}$ is effective,
and so we never have a fiber type contraction of $W_0$. Therefore,
we never have a fiber type contraction of $W$. Thus the minimal
model program ends with an $(W,D)$-minimal model $W\dash W'$
which induces a minimal model $W_0\dash W_0'$ for $(W_0,D|_{W_0})$.
That is, writing $D$ and $L$ for the birational transforms on $W'$
of these divisors,
the pair $(W',D)$ is $\Q$-factorial and klt, $W'$ is terminal, and
$K_{W'}+D\sim_{\Q}cL$ is nef over $T$. (Although we could arrange
for the pair $(W,D)$ to be terminal,
that property may be lost in the course of the minimal model program,
because some components of $D$ may be contracted.)

Because $K_W+D\sim_{\Q}cL$ where $L$ is an integral Weil divisor
and $c>0$, it is a standard property of minimal models that
$H^0(W,O(L))=H^0(W',O(L))$ and $H^0(W_0,O(L))=H^0(W_0',O(L))$
\cite[proof of Theorem 4.5]{DH}. It remains to show that
$H^0(W',O(L))\arrow H^0(W'_0,O(L))$ is surjective.

Since $D$ is big on $W$, it is big on $W'$, and 
so we can write $D=G+H$ for some ample $\Q$-divisor $H$
and effective $\Q$-divisor $G$ on $W'$
such that the pair $(W',G)$ is klt. Then $H^1(W',O_{W'}(L))=0$
follows from Kodaira vanishing, after shrinking the curve
$T$ around 0. (In more detail: we have $L\sim_{\Q}(1/c)(K_{W'}+D)$
where $1/c > 1$. So $L\sim_{\Q} (K_{W'}+G)+(1/c - 1)(K_{W'}+D)+H$,
where $K_{W'}+D$ is nef over $T$ and $H$ is ample over $T$,
which lets us apply Kodaira vanishing on $W'$ over $T$
\cite[Theorem 2.42]{Fujinobook}.)
After shrinking the curve $T$ around 0, the divisor $W_0'$ is linearly
equivalent to zero on $W'$, and so we have $H^1(W',O_{W'}(L-W_0'))=0$.
By de Fernex and Hacon,
we have a short exact sequence
$$0\arrow O_{W'}(L-W_0')\arrow O_{W'}(L)\arrow O_{W_0'}(L)\arrow 0$$
of sheaves on $W'$
\cite[Lemma 4.6]{DH}. Therefore $H^0(W',O(L))\arrow H^0(W_0',O(L))$
is surjective. Equivalently, $H^0(X,O(L))\arrow H^0(X_0,O(L))$
is surjective. \qed

\section{Counterexamples to extension theorems}
\label{counter}

Theorem \ref{class} shows that the divisor class group remains
constant under deformations of a klt Fano variety $X_0$
which is smooth in codimension 2 and $\Q$-factorial in codimension 3.
The first version of this paper asked
whether the Cox ring (of all Weil divisors) deforms
in a flat family under these assumptions. Equivalently, we are asking
whether $H^0(X,O(L))\arrow H^0(X_0,O(L))$ is surjective
for every Weil divisor $L$ on the total space $X$.
In this section, we give
two negative answers to that question. The answer is positive
by de Fernex and Hacon for $\Q$-factorial terminal Fano varieties
\cite[Theorem 1.1, Proposition 6.4]{DH},
and by Theorem \ref{flat} for terminal Fanos which are $\Q$-factorial
in codimension 3. The argument by de Fernex and Hacon actually
shows flatness of Cox rings when we deform
$\Q$-factorial canonical Fano varieties which are smooth
in codimension 2. But the following examples show that the results
mentioned are essentially optimal.

\begin{theorem}
\label{canonical}
There is a canonical Fano 4-fold $X_0$ which is smooth in codimension 3,
a deformation
$X\arrow A^1$ of $X_0$, and a Weil divisor $L$ on $X$
such that $H^0(X,O(L))\arrow H^0(X_0,O(L))$
is not surjective. The fibers $X_t$ for $t\neq 0$ are isomorphic
to $\P^1\times \P^3$.
\end{theorem}

{\bf Proof. }Embed $\P^1\times \P^3$ in projective space $\P^{19}$
by the line bundle $H=O(1,2)$, and let $M$ be a smooth hyperplane section
of $\P^1\times \P^3$. Let $X_0$ be the projective cone over $M$.
Taking $X$ to be a suitable open subset of
the blow-up of the projective cone over $\P^1\times \P^3$ along $M$,
we see that $X_0$ deforms over $A^1$ with
fibers $X_t$ isomorphic to $\P^1\times \P^3$ for $t\neq 0$.

Clearly $X_0$ is smooth in codimension 3. We have $-K_{\P^1\times \P^3}=
O(2,4)=2H$,
and so $-K_M=H$. Therefore the singularity of $X_0$ (the affine cone
over $M$ embedded by $-K_M$) is canonical. Also, $X_0$ is Fano, being
embedded in projective space by an ample line bundle $H$ with
$-K_{X_0}=2H$.

By the Lefschetz hyperplane theorem, the divisor class group of $X_0$ is
$\Cl(X_0)\cong \Pic(M)\cong \Pic(\P^1\times \P^3)\cong \Z^2$.
We have $H^0(X_t,O(-1,2))=0$, but we will show that $H^0(X_0,O(-1,2))
=\oplus_{j\geq 0}H^0(M,O(-1,2)\otimes O(-jH))$ is not zero. It suffices
to show that $H^0(M,O(-1,2))$ is not zero. We have the exact sequence
\begin{eqnarray*}
0  &\arrow  H^0(\P^1\times \P^3,O(i-1,j-2))& \arrow H^0(\P^1\times \P^3,
O(i,j))\arrow H^0(M,O(i,j))\\
  &\arrow  H^1(\P^1\times \P^3,O(i-1,j-2)) &\arrow H^1(\P^1\times \P^3,O(i,j))
\end{eqnarray*}
for any $i,j\in \Z$. By the K\"unneth formula,
$H^1(\P^1\times \P^3,O(-2,0))\cong \C$
and $H^1(\P^1\times \P^3,O(-1,2))=0$. The exact sequence shows
that $H^0(M,O(-1,2))\neq 0$, as we want.
\qed

\begin{theorem}
\label{klt}
There is a $\Q$-factorial klt Fano 4-fold $X_0$ which is smooth
in codimension 3, a deformation $X\arrow A^1$
of $X_0$, and a Weil divisor $L$ on $X$
such that $H^0(X,O(L))\arrow H^0(X_0,O(L))$
is not surjective. The fibers $X_t$ for $t\neq 0$
are isomorphic to the $\P^1$-bundle $P(O\oplus O(3))$ over $\P^3$.
\end{theorem}

{\bf Proof. }The example is similar to Theorem \ref{four},
although the conclusion is different; in that case,
de Fernex-Hacon's extension theorem applies to show
the surjectivity of restriction maps.
First define
$$W=\{ ([x_1,x_2,x_3,x_4,y_2,g],t)\in P(1^4,3,4)\times A^1:
x_1^4+x_2^4+x_3^4+x_4y_2=tg\}.$$
Then $t:W\arrow A^1$ is a flat projective morphism with fibers $W_t$
isomorphic to the weighted projective space $P(1^4,3)$ for $t\neq 0$,
since we can solve the equation for $g$ when $t$ is not zero. Write
$\mu_r$ for the group of $r$th roots of unity. We compute
that $W_0$ has two singular points, the point $[0,0,0,0,1,0]$,
where the singularity is just $A^4/\mu_3$,
and the point $[0,0,0,0,0,1]$. Near the latter
point, $W_0$ is isomorphic to $\{ (x_1,x_2,x_3,x_4,y_2)\in A^5:
x_1^4+x_2^4+x_3^4+x_4y_2=0\}/\mu_4$, where $\mu_4$ acts with
weights $(1,1,1,1,3)$ on these variables. The hypersurface singularity
$x_1^4+x_2^4+x_3^4 +x_4y_2=0$ is terminal, for example by Koll\'ar
\cite[Exercise 67]{Kollarex}, since it has the form
$z_1z_2+f(z_3,\ldots,z_n)=0$ (called a ``$cA$-type singularity'')
and is smooth in codimension 2.
Since a quotient of a klt variety by a finite
group is klt \cite[Proposition 5.20(4)]{KMbook}, $W_0$ is klt.

Also, $W_0$ is $\Q$-factorial by Grothendieck's theorem: a hypersurface
that is smooth in codimension 3 is factorial \cite[Corollaire XI.3.14]{SGA2}.
The same theorem applied to the affine cone over $W_0$ shows that $W_0$
has Picard number 1.

Let $X$ be the blow-up of $W$ along the section $A^1=\{
([0,0,0,0,1,0],t)\}$ of $W\arrow A^1$. Since $W$ is isomorphic to
$(A^4/\mu_3)\times A^1$ near this section, blowing up gets
rid of this singularity. As a result, the morphism
$t:X\arrow A^1$ is smooth outside the single point $([0,0,0,0,0,1],0)$.
Also, the contraction $X\arrow W$ is $K_X$-negative, because
in each fiber of $X\arrow A^1$ we are contracting a $\P^3$
with normal bundle $O(-3)$. Since the fiber $W_t$ is isomorphic to the
weighted projective space $P(1^4,3)$ for $t\neq 0$, $X_t$ is isomorphic
to the $\P^1$-bundle $P(O\oplus O(3))$ over $\P^3$ for $t\neq 0$.
Since the contraction $X_0\arrow W_0$ is an isomorphism near the one
singular point of $X_0$, $X_0$ is smooth in codimension 3,
$\Q$-factorial, and klt.
Finally, $X_0$ has Picard number 2,
since it is a blow-up of $W_0$ at one point.

We now define a flipping contraction $X\arrow Z$ over $A^1$ which contracts
a copy of $\P^1$ in the special fiber $X_0$ to a point. Let
\begin{multline*}
Z=\{([u_{11},u_{12},u_{13},u_{14},u_{21},u_{22},u_{23},u_{24},g],t)\in
P(1^4,4^5)\times A^1: \\
\rank(u_{ij})\leq 1,\, u_{11}^4+u_{12}^4+u_{13}^4+u_{24}=tg\}.
\end{multline*}
Here we view $(u_{ij})$ as a $2\times 4$ matrix to define its rank.
We define the contraction $X\arrow Z$ as the following rational
map $W\dash Z$, which becomes a morphism after blowing up
the above section $A^1$ of $W\arrow A^1$:
$$([x_1,x_2,x_3,x_4,y_2,g],t)\mapsto ([x_1,x_2,x_3,x_4,x_1y_2,x_2y_2,
x_3y_2,x_4y_2,g],t).$$
The resulting contraction $X\arrow Z$ contracts only a copy of
$\P^1$,
the birational transform of the curve $P(3,4)=\{([0,0,0,0,y_2,g],0)\}$
in $W$.

It is straightforward to check that the contraction $X\arrow Z$
is $K_X$-negative; locally, it is a terminal toric flipping contraction,
corresponding to a relation
$e_1+4e_2=e_3+e_4+e_5+e_6$ in $N\cong \Z^5$ in Reid's notation
\cite[Proposition 4.3]{Reid}.
The toric picture shows that the flip $Y\arrow Z$ of $X\arrow Z$
replaces the curve $\P^1$ in the 5-fold $X$ with a copy of $\P^3$.

Since $X_0$ has Picard number 2 and has two nontrivial $K$-negative
contractions (to $W_0$ and $Z_0$), $X_0$ is Fano. Let us show
that $H^0(X,O(L))\arrow H^0(X_0,O(L))$ is not surjective
for some Weil divisor $L$ on $X$. The point is that the flip $X\dash Y$
replaces a curve in the 4-fold $X_0$ by a divisor $\P^3$ in $Y_0$.
Let $A$ be an ample Cartier divisor on $Y$, and also write $A$ for its
birational transform on $X$. Then $H^0(X,O(mA))=H^0(Y,O(mA))$ for each
natural number $m$ because $X$ and $Y$ are isomorphic in codimension one.
Clearly, the images of these groups in $H^0(X_0,O(mA))$
and in $H^0(Y_0,O(mA))$ can be identified.
Define the section ring
$$R(X,A)=\oplus_{m\geq 0} H^0(X,O(mA)).$$
Since $A$ is ample on $Y$,
$H^0(Y,O(mA))$ maps onto $H^0(Y_0,O(mA))$ for $m$ big enough,
and $\Proj(\im(R(Y,A)\arrow R(Y_0,A)))= Y_0$. If the restrictions
$H^0(X,O(mA))\arrow H^0(X_0,O(mA))$ were all surjective, then
the algebra $R(X_0,A)$ would be finitely generated and we would
have $\Proj R(X_0,A)=Y_0$. But that is impossible, since
(when a section ring is finitely generated) the map
$X_0\dash \Proj R(X_0,A)$ is a rational contraction, by Hu and Keel
\cite[Lemma 1.6]{HK}.
(A birational map $X_0\dash Y_0$ is called a rational contraction
if a resolution
$p:Z_0\arrow X_0$ such that $q:Z_0\arrow Y_0$ is a morphism
has the property that every
$p$-exceptional divisor is $q$-exceptional.)
So in fact $H^0(X,O(mA))\arrow H^0(X_0,O(mA))$ is not surjective for some
$m\geq 0$.
\qed


\small \sc DPMMS, Wilberforce Road,
Cambridge CB3 0WB, England

b.totaro@dpmms.cam.ac.uk
\end{document}